\documentclass[12p]{article}

\usepackage{epstopdf}
\usepackage{graphicx}%
\usepackage{float}%
\usepackage{multirow}%
\usepackage{amsmath,amssymb,amsfonts}%
\usepackage{amsthm}%
\usepackage{mathrsfs}%
\usepackage[title]{appendix}%
\usepackage{xcolor}%
\usepackage{textcomp}%
\usepackage{manyfoot}%
\usepackage{booktabs}%
\usepackage{algorithm}%
\usepackage{algorithmicx}%
\usepackage{algpseudocode}%
\usepackage{listings}%
\usepackage{authblk}

\usepackage{hyperref}

\hypersetup{
           breaklinks=true,   
           colorlinks=true,   
           pdfusetitle=true,  
        }

\newtheorem{thm}{Theorem}[section]

\newtheorem{lem}[thm]{Lemma}
\newtheorem{prop}[thm]{Proposition}

\newtheorem{conj}[thm]{Conjecture}
\newtheorem{obs1}[thm]{Observation}
\theoremstyle{definition}

\theoremstyle{remark}
\newtheorem{rem}[thm]{Remark}

\newcommand*{\ditt}{-''-}

\theoremstyle{thmstyleone}%
\theoremstyle{thmstyletwo}%

\theoremstyle{thmstylethree}%
\title{Mutual-visibility sets in Cartesian products of paths and cycles}

\author{Danilo Kor\v ze} 
\affil{Faculty of Electrical Engineering and Computer Science,  University of Maribor,  SI-2000 Maribor, Slovenia }

\author{Aleksander Vesel} 
\affil{Faculty of Natural Sciences and Mathematics, University of Maribor,  SI-2000 Maribor, Slovenia }

\providecommand{\keywords}[1]{\textbf{\textit{Index terms---}} #1}

\begin{document}

\maketitle

\begin{abstract}
For a given graph $G$, the mutual-visibility problem asks for the largest set of vertices $M \subseteq V(G)$ with the property that for any pair of vertices $u,v \in M$ there exists a shortest $u,v$-path of $G$ that does not pass through any other vertex in $M$.

The mutual-visibility problem for Cartesian products of a cycle and a path, as well as for Cartesian products of two cycles, is considered. Optimal solutions are provided for the majority of Cartesian products of a cycle and a path, while for the other family of graphs, the problem is completely solved.
\end{abstract}


\section{Introduction and preliminaries}
In graph theory, a mutual-visibility set is a set of vertices in a graph such that any two vertices in the set are mutually visible, which means that there is a shortest path between them that does not pass through any other vertex in the set. The mutual-visibility number of a graph is the size of a largest mutual-visibility set in the graph. 
The problem of finding the largest mutual-visibility set in the graph can be seen as a relaxation of the 
the general position problem in graphs (for the definition of the problem and the related work see for example \cite{Klavzar} and the references therein).

Mutual-visibility sets have been studied in a variety of contexts, including wireless sensor networks, mobile robot networks, and distributed computing: in wireless sensor networks, mutual-visibility sets can be used to place sensors so that they can communicate with each other without interference; in mobile robot networks, mutual-visibility sets can be used to control robots so that they can avoid collisions; in distributed computing, mutual-visibility sets can be used to design efficient algorithms for problems such as consensus and broadcasting \cite{Aljohani, Bhagat, Cicerone4, DiLuna, Poudel}.

The foregoing problem, with the objective of maximizing the size of the largest mutual-visibility set, was placed on a graph-theoretical footing in 2022 by Di Stefano  \cite{DiS} who shows that this problem is NP-complete and study mutual-visibility sets 
on various classes of graphs, such as block graphs, trees, Cartesian products of paths and cycles,  complete bipartite graphs, and cographs. 

The mutual-visibility in distance-hereditary graphs was studied in \cite{Cicerone3}, 
where it is shown that the mutual-visibility number can be computed in linear time for this class. 
In \cite{Cicerone3}, mutual-visibility sets of triangle-free graphs and Cartesian products were considered. 
In particular, it is shown that computing the mutual-visibility number of a Cartesian product is an intrinsically difficult problem.

Research work from mentioned papers on visibility problems in graphs shows that 
certain modifications to the visibility properties are of interest.  In this respect, \cite{Cicerone2} introduces a variety of new mutual-visibility problems, including the total mutual-visibility problem,  see also \cite{Tian},  the dual mutual-visibility problem, and the outer mutual-visibility problem.

In this paper, we expand the preceding research work on Cartesian products. 
In particular, we study mutual-visibility sets in two families of graphs:
Cartesian products of a cycle and path as well as Cartesian products of two cycles.
In the sequel of this section, we present definitions and results needed for the rest of the paper.
In the following section, we provide some solutions to the mutual-visibility set problem in Cartesian products of a path and 
a cycle. In particular, it is established that the mutual-visibility number for most of the graphs of this class equals  
twofold the length of the cycle. 
In Section 3, we study the mutual-visibility number of Cartesian products of two cycles.
The exact mutual-visibility numbers are provided for all graphs of this class.  Specifically, it is shown that a largest  mutual-visibility set of $C_s \Box C_t$, $s \ge t$,  is of cardinality $3t$  if $t$ is large enough.

Let $G = (V(G),E(G))$ is a graph and $M \subseteq V(G)$.
We say that a shortest $u,v$-path $P$ is {\em $M$-unhindered}, 
 if $P$ does not contain a vertex of $M \setminus \{u, v \}$.
Vertices $u, v \in M$ are 
{\em $M$-visible}  if $G$ admits an $M$-unhindered $u,v$-path.
Moreover,  $M$ is a {\em mutual-visibility set} of $G$ if vertices of $M$ are pairwise $M$-visible. 
The {\em mutual-visibility number}  of $G$ is the size of a largest
mutual-visibility set of $G$ and it is denoted $\mu(G)$.  

The  \emph{Cartesian product} of graphs
$G$ and $H$ is the graph $G \Box H$  with vertex set $V(G) \times V(H)$
and $(x_1,x_2)(y_1,y_2) \in E(G \Box H)$  whenever  $x_1y_1 \in E(G)$ and $x_2=y_2$,  or
 $x_2y_2 \in E(H)$ and $x_1=y_1$.  It is well-known that the Cartesian product is commutative and associative, having the trivial graph as a unit.

The {\em interval} $I(u,v)$ between two vertices $u$ and $v$ of a graph $G$ is the set of vertices on shortest $u,v$-paths. 

For positive integers $n$ and $k$, $n \not = k$, we will use the notation $[k] = \{1, 2, \ldots, k\}$, 
$[k]_0 = \{0, 1, \ldots, k-1\}$ and 
\begin{displaymath}
[k,n] =
        \left \{ \begin{array}{llll}
              \{k, k+1, \ldots, n\}, &  n \ge k   \\
              \{n, n+1, \ldots, k\}, &  k > n    \\
             \end{array}. \right.
\end{displaymath}

For  a positive integer $s$ and $n,k \in [s_0]$,  where $n \not = k$,
  let $$|n - k|_s = \min \{ |n - k|, s - |n - k|\}.$$ 
Moreover, let
\begin{displaymath}
[k,n]_s =
        \left \{ \begin{array}{llll}
              [k,n], &  |n - k|_s <  \frac{s}{2}    \\
              {[s]_0}, &  |n - k|_s =  \frac{s}{2}   \\
              {[s]_0} \setminus [k,n], &  \;  |n - k|_s >  \frac{s}{2}   \\
             \end{array}. \right.
\end{displaymath}

We will assume that $V(P_s)=[s]_0$ for every  $s\ge 2$ and $V(C_t)=[t]_0$ for every  $t\ge 3$.  Thus,  
if $u \in V(X_s \Box X_t)$,   where $X_s$ stands for either  $P_s$ or $C_s$,  while
$u_x \in [s]_0$ and  $u_y \in [t]_0$, we write $u=(u_x,u_y)$. 
If $i \in [t]_0$ (resp.  $j \in [s]_0$), 
then the subgraph of $C_s \Box C_t$  induced by $V (C_s) \times \{i\}$ (resp. $ \{j\} \times V (C_t)   $)
is isomorphic to $C_s$ (resp.  $C_t$); it is called a {\em$C_s$-fiber}  
(resp.  {\em$C_t$-fiber}) and  denoted $C_s^i$ (resp. $^j\!C_s$).

A $P_s$-fiber and $C_t$-fiber in $P_s \Box C_t$ denoted by $P_s^i$ and $^j\!C_s$, respectively,  are defined analogously.

The following result can be easily confirmed.  
\begin{obs1} \label{o1}
Let $t$ and $s$ be nonnegative integers.

1.  If $u,v \in V(P_s \Box P_t)$, then $I(u,v) = [u_x,v_x]  \times  [u_y,v_y]$.

2.  If $u,v \in V(P_s \Box C_t)$, then $I(u,v) = [u_x,v_x]  \times  [u_y, v_y]_t$.

3.  If $u,v \in V(C_s \Box C_t)$, then $I(u,v) = [u_x,v_x]_s  \times  [u_y,v_y]_t$.
\end{obs1}

Note also that for  $u,v \in V(P_s \Box P_t)$, $u',v' \in V(P_s \Box C_t)$  and  $u'',v'' \in V(C_s \Box C_t)$ 
we have $d(u,v) = |u_x - v_x| +  |u_y - v_y|$, $d(u',v') = |u_x' - v_x'| +  |u_y' - v_y'|_t$ 
and $d(u'',v'') = |u_x'' - v_x''|_s +  |u_y'' - v_y''|_t$. 

\begin{lem} \label{l1}
Let $M \subseteq V(P_s \Box P_t)$. If $M$ obeys the following conditions
 
 (i) $|M \cap P_s^i| \le 2 $ and $ |M \cap ^jP_t| \le 2$ for every  $i \in [t]_0$ and $j \in [s]_0$,
  

(ii) if  $M \cap P_s^i = \{ u, v\}$ and $M \cap P_s^j = \{ w,  z\}$ for some $i,j \in [t]_0$,  then \\  \hphantom{ \ \ \ \ \ \ \ }  
$[u_x,v_x]  \cap [w_x,z_x] \not = \emptyset$, 

(iii) if  $M \cap ^i\!P_t = \{ u, v\}$ and $M \cap ^j\!P_t = \{ w,  z\}$ for some $i,j \in [s]_0$,  then \\  \hphantom{ \ \ \ \ \ \ \ }  
$[u_y,v_y]  \cap [w_y,z_y] \not = \emptyset$, 
  
 (iv) $d(u,v) \ge 3$ for every $u,v \in M$, 
 
\noindent
then  $M$ is a  mutual-visibility set of $P_s \Box P_t$.
\end{lem}

\begin{proof}
We will show that there exists an $M$-unhindered  $u,v$-path for every  $u,v \in M$. 
Note that for $u_x = v_x$ or $u_y = v_y$, by condition (i) of the lemma, the assertion trivially holds. 
Let assume w.l.o.g. that $u_y < v_y$ and let $M_{u,v} := (M \setminus \{u, v \}) \cap [u_x,v_c] \times [u_y, v_y]$. 

The proof is by induction on the cardinality of $M_{u,v}$.  
If $|M_{u,v}|=1$,  then $u_x \not = v_x$  by condition (i) of the lemma.  Thus,  it is not difficult to construct a shortest $u,v$-path that does not contain the vertex of   $M_{u,v}$. 
 
Suppose now that the lemma does not hold and let $u$ and $v$ be vertices of $M$ which are not mutually visible in $P_s \Box P_t$.  Moreover, let $u$ and $v$ be chosen in the way that $M_{u,v}$ is the set with the smallest cardinality 
such that $u$ and $v$ possess this property.  Suppose first that $u_x < v_x$. 
Note that for every $w \in M_{u,v}$,  by the minimality of $M_{u,v}$, there exists an $M$-unhindered  $w,v$-path, say $P$.   
Note that,  by $w_x \le v_x$ and $w_y \le v_y$, it holds that $P$ admits either 
the vertex $w'=(w_x+1,w_y)$ or the vertex $w''=(w_x,w_y+1)$.  

If there exists $w \in M_{u,v}$ such that $M_{u,w} = \emptyset$,  $u_x \not = w_x$ and  $u_y \not = w_y$,   let  $w_x - u_x = k$ and $w_y- u_y = \ell$. We will show that we can always find an $M$-unhindered  $u,w'$- and $u,w''$-path. 
Since $M$ is composed of  vertices that are pairwise at a distance of at least three,  
$P'=u, (u_x+1,u_y), \ldots,(u_x+k,u_y),(u_x+k,u_y+1),\ldots,
(u_x+k,u_y+\ell-1),  (u_x+k+1,u_y+\ell-1),  (u_x+k+1,u_y+\ell) =w'$ is an  $M$-unhindered $u,w'$-path,  while
$P''=u, (u_x+1,u_y), \ldots,(u_x+k-1,u_y),(u_x+k-1,u_y+1),\ldots,
(u_x+k-1,u_y+\ell+1),  (u_x+k,u_y+\ell+1) =w''$ is an $M$-unhindered  $u,w''$-path.

Otherwise  (i.e. , $M$ does not admit $w$ from the previous paragraph),
suppose w.l.o.g.  that there exists $z^1 \in M_{u,v}$ such that $z^1_y = u_y$.  
Let us consider the sequence of vertices $z^1,  z^2, \ldots,  z^{a}$, $a \ge 1$, such that 
$z^{i+1}_y = z^{i}_y + 1$ and $z^{i+1}_x > z^{i}_x$,  
 where $a$ is the largest integer that allows this sequence.  
 Note that  from condition (ii) of the lemma it follows that 
 $[z^1_x,  z^a_x] \times [z^1_y,  z^a_y] \cap M \subseteq  \{z^{1}, \ldots, z^{a} \}$.  
 
If there exists  $w \in M_{u,v}$ and $i \in [ a]$
such that $w_y > z^{a}_y$,  $ z^{i+1}_x \ge w_x \ge z^{i}_x$ and  $M_{z^{i},w} = \emptyset$,
then by the above discussion  we have  $M_{u,w} = \{z^{1}, \ldots, z^{i} \}$.
Since $M$ is composed of vertices that are pairwise at a distance of at least three,  
 we can now construct, similarly as above, an $M$-unhindered  $u,w'$- and $u,w''$-path.
 
Otherwise (i.e.,  $M$ does not admit $w$ from the previous paragraph),   we clearly have $z^{a} = v$.  But since it is not difficult to construct an $M$-unhindered  $u,v$-path for this case we obtain a contradiction. 
 
Note that we showed above that we can always construct an $M$-unhindered  $u,w'$- and $u,w''$-path.
If $w' \in V(P)$,  then, $V(P) \setminus  \{ w \}$ is an $M$-unhindered   $w',v$-path.  Moreover,  the graph induced by 
$V(P) \cup V(P') \setminus  \{ w \}$ is an $M$-unhindered $u,v$-path,  say $Q$.  
 Since $Q$ is of length $w_x + 1 - u_x + v_x - w_x - 1 + w_y-u_y+v_y-w_y = v_x - u_x + v_y- u_y$ it holds that $Q$ is an $M$-unhindered $u,v$-path and we obtain a contradiction.  If $w'' \in V(P)$,  then we can analogously construct  
an $M$-unhindered $u,v$-path that yields a contradiction. 

Since the proof for $u_x > v_x$ is analogous, the proof is complete.
\end{proof}

\section{Cartesian products of  a cycle and a path}

By $\mu(P_s ) = 2$ and  $\mu(C_t ) = 3$ for every $t>s\ge 2$,   the following result follows from  \cite[Lemma 2.3]{DiS}.
\begin{prop} \label{p2}
If $s \ge 2$ and $t \ge 3$,  then $\mu(P_s \Box C_t) \le min\{3s,  2t \}$.
\end{prop}

\begin{figure}[H]
     \centering
	\includegraphics[width=7cm]{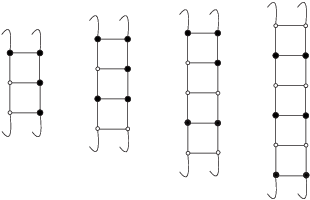}
		\bigskip
	\caption{Mutual-visibility sets of $P_{2} \Box  C_{3}$, $P_{2} \Box  C_{4}$, $P_{2} \Box  C_{5}$ and $P_{2} \Box  C_{6}$}
\label{P2}
\end{figure}

If the path in a product is of length 2, we obtain the following result.

\begin{prop} \label{p3}
If $t \ge 3$, then 
\begin{displaymath}
\mu(P_2 \Box C_t) =
        \left \{ \begin{array}{llll}
              4,  &  t = 3   \\
              5,  &  t = 4   \\
              5,  &  t = 5   \\
              6,  &  t \ge 6   \\
             \end{array}. \right.
\end{displaymath}
\end{prop}

\begin{proof}
By Proposition \ref{p2}, we have  $\mu(P_2 \Box C_t) \le 6$. 
For $t \le 6$, the result follows from the constructions depicted in Fig. \ref{P2} and straightforward case analysis.  
Since it is not difficult to provide construction with 6 vertices for every $t \ge 7$, the proof is complete.
\end{proof}

The next observation can be easily confirmed.  

\begin{obs1} \label{o2}
Let $t$ and $s > k$ be nonnegative integers.  
If $M$ is a mutual-visibility set of $P_k \Box C_t$, then $M$ is a mutual-visibility set of $P_s \Box C_t$. 
\end{obs1}

 The following  proposition shows that the upper bound  on the mutual-visibility number of $P_s \Box C_t$
is  achieved when $t$ and $s$ are big enough. 

\begin{thm} \label{t1}
If  $s + 1 \ge t \ge 6 $,  then $\mu(P_s \Box C_t) = 2t$.
\end{thm}

\begin{figure}[H]
     \centering
	\includegraphics[width=12cm]{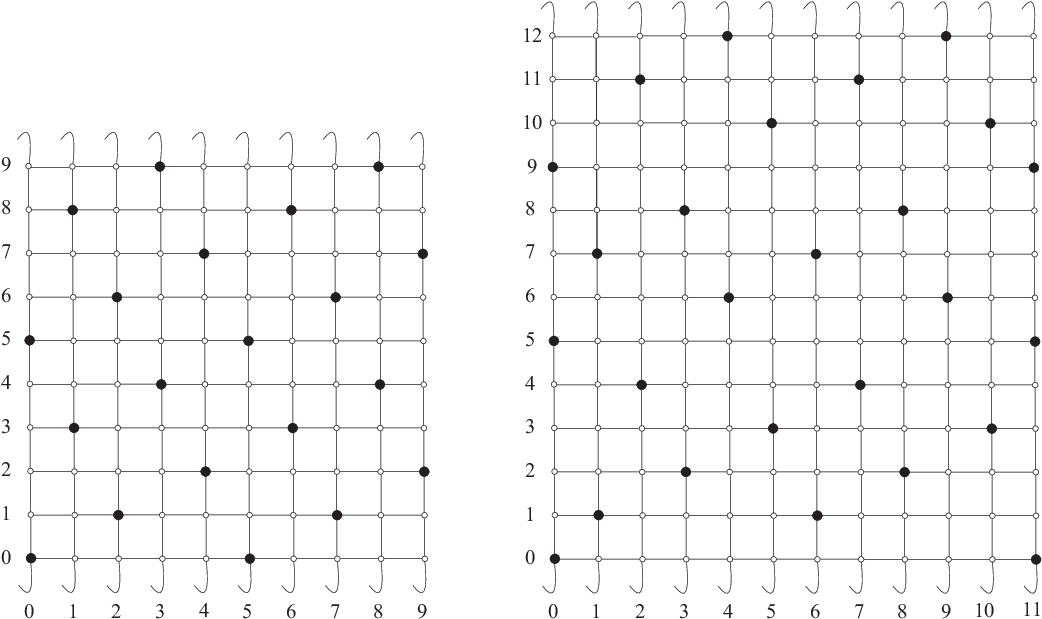}
		\bigskip
	\caption{Mutual-visibility set of $P_{10} \Box  C_{10}$ (left) and   $P_{12} \Box  C_{13}$  (right)}
\label{P10C10}
\end{figure}

\begin{proof}
By Observation \ref{o1},  for every $u,v \in V(P_{s} \Box C_{t})$  we have  $I(u,v) = [u_x,v_x] \times [u_y,v_y]_{t}$. Thus,  the graph induced by $[u_x,v_x]  \times  [u_y,v_y]_{t}$  is isomorphic to a subgraph of $P_{s}  \Box P_{\lceil {t+1 \over 2}\rceil}$.

Let $t \ge 13$ is odd and $M' =  \{((2i + j (2 \lfloor {t-3 \over 4}\rfloor + 1)) \, {\rm mod} \, (t-3), i) \, | \, i \in [t-3]_0,  j\in [2]_0 \}$.  This set for $t=13$ is depicted on the left-hand side of Fig. \ref{P10C10}. 
Obviously,  $M'$ is a subset of $P_{t-3} \Box C_{t-3}$ of cardinality $2(t-3)$. 
By confirming the conditions of Lemma \ref{l1} we will show that $M'$ is a mutual-visibility set of $P_{t-3} \Box C_{t-3}$.

Note first that $M' =  \{(2i,  (i  + j {t-3 \over 2}) \, {\rm mod} \, (t-3)) \, | \, i \in [{t-3 \over 2}]_0,  j\in [2]_0 \}
\cup \{(2i+1,  (i  +\lceil {t-3 \over 4}\rceil +  j {t-3 \over 2}) \, {\rm mod} \, (t-3)) \, | \, i \in [{t-3 \over 2}]_0,  j\in [2]_0 \}$. 
Thus, it readily follows that 
$|M' \cap P_{t-3}^k| = 2 $ and $ |M' \cap ^{\ell}\!C_{t-3}| = 2$ for every  $k, \ell \in [t-3]_0$ which in turn satisfies
condition (i) of Lemma \ref{l1}.
Moreover, it is not difficult to see that $d(u,v) \ge 3$ for every $u,v \in M'$.  Thus,  we are left with conditions (ii) and (iii) of Lemma \ref{l1}.
By the above discussion, for every $u,v \in V(P_{t-3} \Box C_{t-3})$ we have $I(u,v) = [u_x,v_x] \times [u_y,v_y]_{t-3}$.  It follows that the graph induced by $[u_x,v_x]  \times  [u_y,v_y]_{t-3}$  is isomorphic to a subgraph of 
$P_{t-3}  \Box P_{\lceil {t-2 \over 2}\rceil}$.
Since for every $u,v \in M' \cap P_{t-3}^k$, $k \in [t-3]_0$, it holds that $d(u,v)= 2 \lceil {t-3 \over 4}\rceil + 1$,
it follows that for every $u,v \in M' \cap P_{t-3}^k$ and every $w,z \in  M' \cap P_{t-3}^{k'}$ we have 
$\min \{u_x, v_x \} < \max \{z_x, w_x \}$.
Moreover,  since
for every $u,v \in M' \cap ^\ell\!C_{t-3}$, $\ell \in [t-3]_0$ it holds that $d(u,v)=  \lceil {t-3 \over 2}\rceil $,
it follows that for every $u,v \in M' \cap ^{\ell}\!C_{t-3}$ and every $w,z \in  M' \cap ^{\ell'}\!C_{t-3}$ we have 
$\min \{u_y, v_y \} < \max \{z_y, w_y \}$.
 Thus, conditions (ii) and (iii) of Lemma \ref{l1} are also fulfilled. 

Let $M := B \cup M_1 \cup M_2 \cup M_3$,  where the sets $B$,  $M_1$,  $M_2$,  and
$M_3$ are defined as follows:

$B =\{(0,0),(0,\lceil {t \over 3}\rceil), (0, \lceil {2 t \over 3}\rceil),(t-2,0),
(t-2,\lceil {t \over 3}\rceil), (t-2, \lceil {2 t \over 3}\rceil) \}$,

$M_1 = \{((2i + j (2 \lfloor {t-3 \over 4}\rfloor + 1)) \, {\rm mod} \, (t-3)+1, i+1) \, | \, i \in [\lceil {t \over 3}\rceil-1]_0,  j\in [2]_0 \}$,

$M_2 = \{((2i + j (2 \lfloor {t-3 \over 4}\rfloor + 1)) \, {\rm mod} \, (t-3)+1, i+2) \, | \, i \in [\lceil {t \over 3}\rceil-1,
 \lceil {2 t \over 3}\rceil-3],  j\in [2]_0 \}$,

$M_3 = \{((2i + j (2 \lfloor {t-3 \over 4}\rfloor + 1)) \, {\rm mod} \, (t-3)+1, i+3) \, | \, i \in [\lceil {2 t \over 3}\rceil-2,
t-4],  j\in [2]_0 \}$,

Note that $M$ is a subset of $P_{t-1} \Box C_{t}$ and remind that $M'$ is a mutual-visibility set
of $P_{t-3} \Box C_{t-3}$. 
We can see that $M$ can be obtained from $M'$ by inserting in $P_{t-3} \Box C_{t-3}$
three $P_ {t-1}$-fibers and two $C_ {t}$-fibers that contains vertices of $B$ 
(for $t =13$ see the right-hand side of Fig. \ref{P10C10}).  
Thus,  $M \setminus B$ is clearly a mutual-visibility set of $P_{t-1} \Box C_{t}$.
 Moreover,  the vertices of $B\cup M$ do not violate the conditions of Lemma \ref{l1} with the exception of 
condition (iv), i.e.,  a vertex of $B $ may be at distance two from a vertex of $M \setminus B$. 
However, since a vertex  $u \in B$ belongs either to   $^{0}C_{t}$ or $^{t-1}\!C_{t}$ it is not difficult 
to obtain an $M$-unhindered  $u,v$-path for every $v \in M$.  
Suppose for example that $u  \in \, ^{0}\!C_{t}$,  $v, w \in M \setminus B$,  $d(u,w)=2$ and $u_y < w_y \le  v_y$.
If $P$ denote an $M$-unhindered $w,v$-path,  then either $w'=(w_x+1, w_y)$ or $w''=(w_x, w_y+1)$ belongs to $P$.
Since it is trivial to obtain an $M$-unhindered $u,w'$-path and $u,w''$-path,  we can also obtain 
an $M$-unhindered  $u,v$-path.
 
It follows that $M$ is a mutual-visibility set of $P_{t-1} \Box C_{t}$. 

\begin{figure}[H]
     \centering
	\includegraphics[width=12cm]{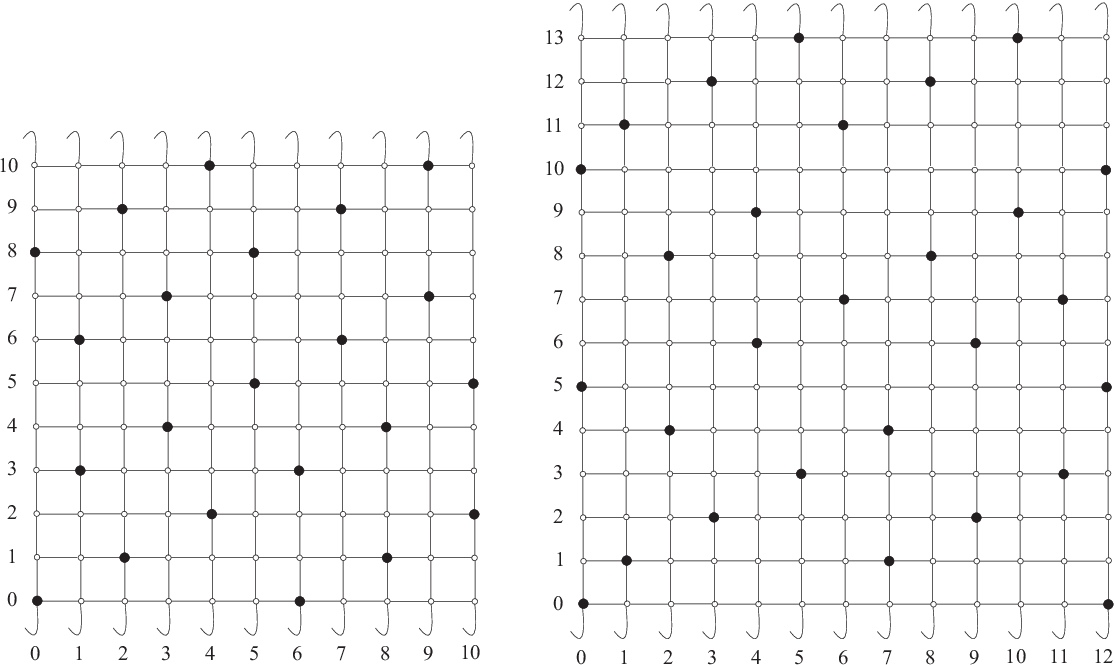}
		\bigskip
	\caption{Mutual-visibility set of (left) $P_{11} \Box  C_{11}$ and  $P_{13} \Box  C_{14}$ (right)}
\label{P11C12}
\end{figure}

\begin{table}[H]  
\begin{center}
\caption{Mutual-visibility numbers of $P_s \Box C_t$,  where $t \le 17 $}  \label{tab1}
\bgroup
\def\arraystretch{1.5}
\begin{tabular}{|c||c|c|c|c|c|c|c|c|c|c|}
  \hline
  $t$ \textbackslash $s$ & $3$ & $4$ & $5$ & $6$ & $7$  &  $8$  &  $9$ &  $10$ &  $11$ &  $12$ \\ 
  \hline
 \hline
  $3$ & 6  & \ditt  & \ditt  & \ditt &  \ditt & \ditt &  \ditt &  \ditt  & \ditt & \ditt \\
\hline
  $4$ & $7$  & 8         & \ditt  & \ditt&  \ditt & \ditt &  \ditt &  \ditt & \ditt & \ditt\\
\hline
  $5$ &  7 &  $9$   & $10$  & \ditt &  \ditt & \ditt &  \ditt &  \ditt & \ditt & \ditt\\
\hline
 $6$ & $8$  & $10$ & $12$  & \ditt &  \ditt & \ditt &  \ditt &  \ditt & \ditt & \ditt\\
\hline
  $7$ & $8$  & $10$  & $12$  & $14$ &  \ditt & \ditt &  \ditt &  \ditt & \ditt & \ditt\\
\hline
  $8$ & $9$  & $12$  & $14$  & $16$ &  \ditt & \ditt &  \ditt &  \ditt & \ditt & \ditt\\
\hline
  $9$ & $9$  & $12$  & $14$  & $16$ &  $17$ & $18$ &  \ditt &  \ditt & \ditt & \ditt\\
    \hline
  $10$ & $9$  & $12$  & $15$  & $18$ &  $19$ & $20$ &  \ditt &  \ditt & \ditt & \ditt\\
    \hline
  $11$ & $9$  & $12$  & $15$  & $18$ &  $19$ & $22$ &  \ditt &  \ditt & \ditt & \ditt\\
    \hline
  $12$ & $9$  & $12$  & $15$  & $18$ &  $21$ & $24$ &  \ditt &  \ditt & \ditt & \ditt\\
    \hline
  $13$ & $9$  & $12$  & $15$  & $18$ &  $21$ & $24$ &  $26$ &  \ditt & \ditt & \ditt\\
    \hline
  $14$ & $9$  & $12$  & $15$  & $18$ &  $21$ & $24$ &  $27$ &  $28$ & \ditt & \ditt\\
    \hline
  $15$ & $9$  & $12$  & $15$  & $18$ &  $21$ & $24$ &  $27$ &  $30$ & \ditt& \ditt \\
    \hline
  $16$ & $9$  & $12$  & $15$  & $18$ &  $21$ & $24$ &  $27$ & $30$ & $32$ & \ditt \\
    \hline
  $17$ & $9$  & $12$  & $15$  & $18$ &  $21$ & $24$ &  $27$ &  $30$ &  $33$ & $34$ \\
    \hline
\end{tabular}
\egroup
\end{center}
\end{table}

Let $t \ge 14$ is even and 
$M' =  \{((2i + j (t-8)) \, {\rm mod} \, (t-3), i) \, | \, i \in [t-3]_0,  j\in [2]_0 \}$. 
This set for $t=14$ is depicted on the left-hand side of Fig. \ref{P11C12}. 
Obviously,  $M'$ is a subset of $P_{t-3} \Box C_{t-3}$ of cardinality $2(t-3)$. 
We will show that $M'$ is a mutual-visibility set of $P_{t-3} \Box C_{t-3}$ by confirming the conditions 
of Lemma \ref{l1}.

We can see that  $M' =  \{(i,  (i {t-2 \over 2} +j {t+2 \over 2} ) \, {\rm mod} \, (t-3)) \, | \, i \in [t-3]_0,  j\in [2]_0 \}$.
Thus, it readily follows that 
$|M' \cap P_{t-3}^k| = 2 $ and $ |M' \cap ^{\ell}\!C_{t-3}| = 2$ for every  $k, \ell \in [t-3]_0$ which in turn satisfies
condition (i) of Lemma \ref{l1}.
Moreover, it is not difficult to see that $d(u,v) \ge 3$ for every $u,v \in M'$.  In order to see that  
conditions (ii) and (iii) of Lemma \ref{l1} are fulfilled, 
 note that for every $u,v \in M' \cap P_{t-3}^k$, $k \in [t-3]_0$, it holds that $d(u,v)= t-8$, 
 while for for every $u,v \in M' \cap ^\ell\!C_{t-3}$, $\ell \in [t-3]_0$,  it holds that 
 $d(u,v)=  \lceil {t+2 \over 2}\rceil $. 
 Remind that the graph induced by $[u_x,v_x]  \times  [u_y,v_y]_{t-3}$  is isomorphic to a subgraph of 
$P_{t-3}  \Box P_{\lceil {t-2 \over 2}\rceil}$.
Analogously as for an even $t$ we can then conclude that conditions (ii) and (iii) of Lemma \ref{l1} are also fulfilled. 
It follows that $M'$ is a mutual-visibility set of $P_{t-3} \Box C_{t-3}$.

Let $M := B \cup M_1 \cup M_2 \cup M_3$,  where the sets $B$,  $M_1$,  $M_2$,  and
$M_3$ are defined as follows:

$B =\{(0,0),(0,\lceil {t \over 3}\rceil), (0, \lceil {2 t \over 3}\rceil),(t-2,0),
(t-2,\lceil {t \over 3}\rceil), (t-2,\lceil {2 t \over 3}\rceil) \}$,

$M_1 =  \{((2i + j (t-8)) \, {\rm mod} \, (t-3) + 1,  i + 1) \, | \,  i \in [\lceil {t \over 3}\rceil-1]_0,  j\in [2]_0 \}$.

$M_2 =  \{((2i + j (t-8)) \, {\rm mod} \, (t-3) + 1,  i + 2) \, | \,  i \in [\lceil {t \over 3}\rceil-1,
 \lceil {2 t \over 3}\rceil-3],    j\in [2]_0 \}$,

$M_3 =  \{((2i + j (t-8)) \, {\rm mod} \, (t-3) + 1,  i + 3) \, | \,  i \in [\lceil {2 t \over 3}\rceil-2,
t-4],   j\in [2]_0 \}$,

Again,  we can see that $M$ can be obtained from $M'$ by inserting in $P_{t-3} \Box C_{t-3}$
three $P_ {t-1}$-fibers and two $C_ {t}$-fibers that contains vertices of $B$ and 
analogously as for an even $t$ we can conclude that $M$ is a mutual-visibility set of $P_{t-1} \Box C_{t}$ if $t$ is even
(for $t =14$ see the right-hand side of Fig. \ref{P11C12}). 

From the above discussion, it follows that there exists a  mutual-visibility set with $2t$ vertices of $P_{t-1} \Box C_{t}$
for every $t \ge 13$.  Moreover,   since we found a mutual-visibility set of this size for every $6 \le t \le 12$
(these sets  are available on the web page 
\href{https://omr.fnm.um.si/wp-content/uploads/2023/08/MVS-constructions-for-paths-and-cycles.txt}{external document}), we have
$\mu(P_{t-1} \Box C_t) = 2t$ for every $t \ge 13$.  
Proposition \ref{p2} and Observation \ref{o2} now yield the assertion.
\end{proof}

\begin{rem}  
Mutual-visibility numbers of $P_s \Box C_t$,  for $t \le 12 $,  obtained by a computer,  are presented in 
Table   \ref{tab1}.  Note that the last shown mutual-visibility number in every row equals $2t$.
\end{rem}

It can be seen that for every $12 \le t \le 17$ we have   $\mu(P_s \Box C_t) = min\{3s,  2t \}$. 
This observation together with Proposition \ref{p3} leads us to state the following conjecture.

\begin{conj}
If $t \ge 12$,  then $\mu(P_s \Box C_t) = min\{3s,  2t \}$.
\end{conj}

\section{Cartesian products of two cycles}

The following result is given in \cite{DiS}.
\begin{prop} \label{p4}
If  $s \ge t \ge 3$,  then $\mu(C_s \Box C_t) \le 3t$.
\end{prop}

The next results show that the upper bound is achieved for most products of the form $C_t \Box C_t$.
By Proposition \ref{p3},  it holds $\mu(C_t \Box C_t) \le 3t$. Thus,  we have to show that a mutual-visibility set of cardinality  $3t$ exists for every $t $.

\begin{prop} \label{mod3}
Let $t\equiv 3$ (mod $6$).  If $t \ge 15$,  then $\mu(C_t \Box C_t) = 3t$.
\end{prop}

\begin{figure}[H]
     \centering
	\includegraphics[width=12cm]{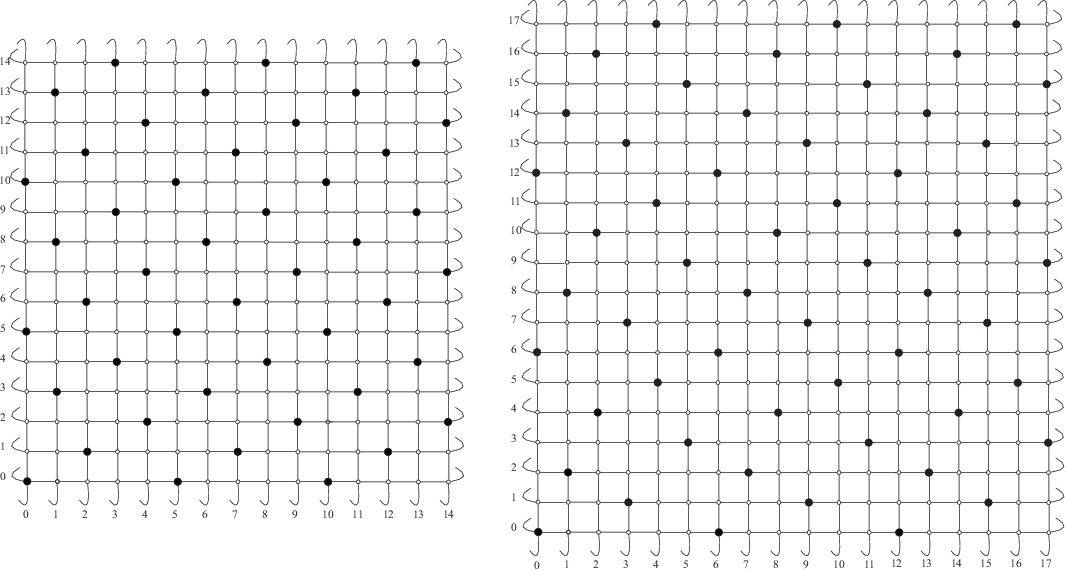}
		\bigskip
	\caption{Mutual-visibility set of $C_{15} \Box  C_{15}$ (left) and  $C_{18} \Box  C_{18}$ (right)}
\label{C15C15}
\end{figure}

\begin{proof}
We will show that a mutual-visibility set of cardinality  $3t$
exists for every $t \equiv 3$ (mod $6)$, $t \ge 15$.  
Let $t=3k$, where $k\ge 5$ is odd.
Note that for every $u,v \in V(C_{3k} \Box C_{3k})$  we have 
$I(u,v) = [u_x,v_x]_{3k} \times [u_y,v_y]_{3k}$. Moreover,  the graph induced by
$[u_x,v_x]_{3k} \times [u_y,v_y]_{3k}$  is isomorphic to a subgraph of 
$P_{  { 3k+1 \over 2  }}  \Box P_{{3k+1 \over 2}}$.

Let $M = \{((2i + jk) \, {\rm mod} \, 3k, i) \, | \, i \in [3k]_0,  j\in [3]_0 \}$.
For $k=5$, the construction is depicted on the left-hand side of Fig.  \ref{C15C15}. 
We can see that  $M = \{(i, ({k +1 \over 2} i + jk) \, {\rm mod} \, 3k) \, | \, i \in [3k]_0,  j\in [3]_0 \}$.

Since the vertices of $M \cap C_{3k}^i $ as well as of  $M \cap ^{j}\!C_{3k}$ are pairwise at distance $k$,  it readilly follows that 
$|M \cap C_{3k}^i \cap [z_x,w_x]_{3k} \times [z_y,w_y]_{3k}| \le 2 $ and 
$|M \cap ^{j}\!C_{3k}\cap [z_x,w_x]_{3k} \times [z_y,w_y]_{3k}| \le 2$ for every  $i,j \in [3k]_0$ and every 
$z,w \in V(C_{ 3k}  \Box C_{3k})$, which in turn satisfies condition (i) of Lemma \ref{l1}.
Moreover, since 
$d(z_x,w_x)  \le {3k - 1 \over 2}$ and $d(z_y,w_y)  \le {3k - 1 \over 2}$ for every 
$z,w \in V(C_{ 3k}  \Box C_{3k})$, 
conditions (ii) and (iii) of Lemma \ref{l1} are also fulfilled. 

We can also see that  $d(u,v) \ge 3$ for every $u,v \in M$,  which confirms condition (iv) of Lemma \ref{l1}.
It follows that $M $ is a mutual-visibility set of cardinality  $9k$ in $C_{ 3k}  \Box C_{3k}$ if $k\ge 5$ is odd.
\end{proof}
 
 \begin{prop} \label{mod0}
Let $t\equiv 0$ (mod $6$).  If $t \ge 18$,  then $\mu(C_t \Box C_t) = 3t$.
\end{prop}
 
\begin{proof}
We will  show that a mutual-visibility set of cardinality  $3t$
exists for every $t \equiv 0$ (mod $3)$, $t \ge 18$.   Let $t=3k$, where $k\ge 6$ is even.
Note that for every $u,v \in V(C_{3k} \Box C_{3k})$  we have 
$I(u,v) = [u_x,v_x]_{3k} \times [u_y,v_y]_{3k}$. Moreover,  the graph induced by
$[u_x,v_x]_{3k} \times [u_y,v_y]_{3k}$  is isomorphic to a subgraph of 
$P_{  { 3k+2 \over 2  }}  \Box P_{{3k+2 \over 2}}$.

Let  $M =  \{((i  + j k, 2i + \ell k) \, | \, i \in [{k  \over 2}]_0,  j, \ell \in [3]_0 \}
\cup \{(k-i-1+jk,  (2i  + 3 + \ell k ) \, {\rm mod} \, 3k) \, | \, i \in [{k \over 2}]_0,  j, \ell \in [3]_0 \}$. 

For $k=6$, the construction is depicted on the right-hand side of Fig.  \ref{C15C15}. 
We can see that the vertices of $M \cap C_{3k}^i $ as well as of  $M \cap ^{j}\!C_{3k}$ are pairwise at distance $k$.
Thus,  using the same arguments as in the proof of Proposition \ref{mod3}  we can see that $M$ satisfies conditions (i), (ii), and (iii) of Lemma \ref{l1}.

Since it is not difficult to confirm that $d(u,v) \ge 3$ for every $u,v \in M$,  
it follows that $M $ is a mutual-visibility set of cardinality  $9k$ in $C_{ 3k}  \Box C_{3k}$ if $k\ge 6$ is even.
\end{proof}

\begin{figure}[H]
     \centering
	\includegraphics[width=12cm]{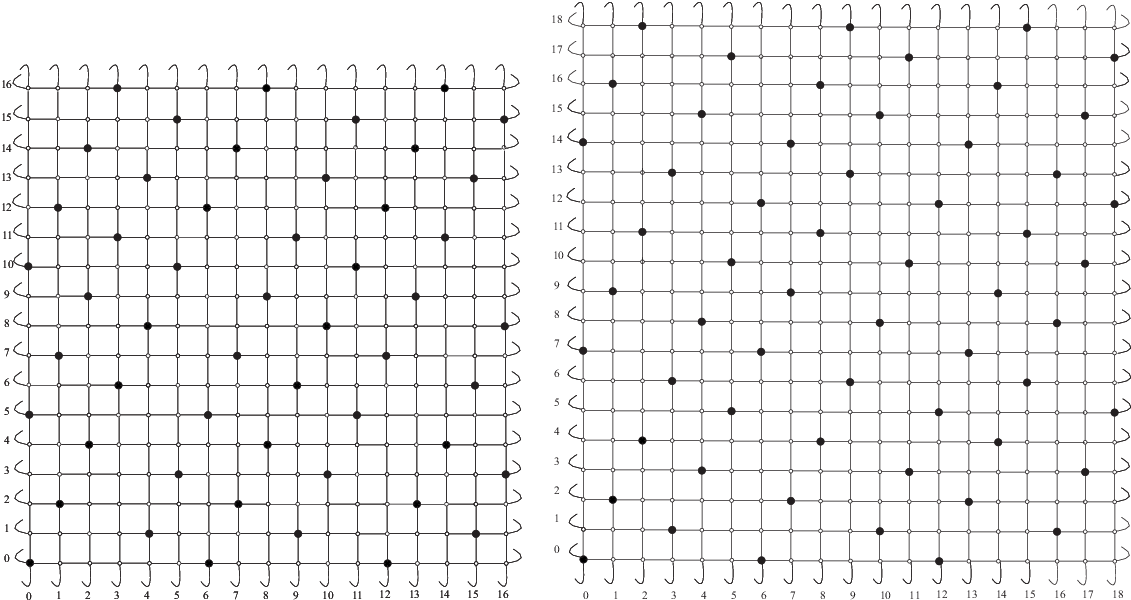}
		\bigskip
	\caption{Mutual-visibility set of $C_{17} \Box  C_{17}$ (left) and  $C_{19} \Box  C_{19}$ (right)}
\label{C17C17}
\end{figure}

\begin{prop} \label{mod5}
Let $t\equiv 5 $ (mod $6$).  If $t \ge 17$,  then $\mu(C_t \Box C_t) = 3s$.
\end{prop}

\begin{proof}
We will show that a mutual-visibility set of cardinality  $3t$ exists for every $t \equiv 5$ (mod $6)$, $t \ge 17$.  

Let us of an odd $k$ set $t=3k+2$.
Note that for every $u,v \in V(C_{3k+2} \Box C_{3k+2})$  we have 
$I(u,v) = [u_x,v_x]_{3k+2} \times [u_y,v_y]_{3k+2}$. Moreover,  the graph induced by
$[u_x,v_x]_{3k+2} \times [u_y,v_y]_{3k+2}$  is isomorphic to a subgraph of 
$P_{ {3k+3 \over 2}}  \Box P_{3k+3 \over 2}$.

Let 
$M =  \{(i  + j (k+1) \, {\rm mod} \, (3k+2)),  2i ) \, | \, i \in [{3k+3  \over 2}]_0,  j\in [3]_0 \}
\cup \{(( i  + {3k+3 \over 2} +  j (k+1) \, {\rm mod} \, (3k+2)), 2i+1) \, | \, i \in [{3k+1 \over 2}]_0,  j\in [3]_0 \}$. 
For $k=5$ the construction is depicted on the left-hand side of Fig.  \ref{C17C17}. 
We can see that
$M = \{(i, (2i + jk) \, {\rm mod} \, (3k+2)) \, | \, i \in [3k+2]_0,  j\in [3]_0 \} $.

Note that the distance between vertices of $M \cap C_{3k+2}^i $ is either $k$ or  $k+1$, while the distance between vertices of  $M \cap ^{j}\!C_{3k+2}$ is either $k$ or  $k+2$.
Thus,  using the same arguments as in the proof of Proposition \ref{mod3}  we can see that $M$ satisfies conditions (i), (ii), and (iii) of Lemma \ref{l1}.

Since we can confirm that $d(u,v) \ge 3$ for every $u,v \in M$,  
It follows that $M $ is a mutual-visibility set of cardinality  $9k+6$ in $C_{ 3k+2}  \Box C_{3k+2}$ if $k\ge 5$ is odd.
\end{proof}

\begin{prop} \label{p4}
Let $t\equiv 1 $ (mod $6$).  If $t \ge 19$,  then $\mu(C_t \Box C_t) = 3t$.
\end{prop}

\begin{proof}
We will show that a mutual-visibility set of cardinality  $3t$ exists for every $t \equiv 1$ (mod $6)$,   $t \ge 19$. 

Let us of an even $k$ set $t=3k+1$.
Note that for every $u,v \in V(C_{3k+1} \Box C_{3k+1})$  we have 
$I(u,v) = [u_x,v_x]_{3k+1} \times [u_y,v_y]_{3k+1}$. Moreover,  the graph induced by
$[u_x,v_x]_{3k+1} \times [u_y,v_y]_{3k+1}$  is isomorphic to a subgraph of 
$P_{ {3k+2 \over 2}}  \Box P_{3k+2 \over 2}$.

Let 
$M =  \{((i  + j k) \, {\rm mod} \, (3k+1)),2i) \, | \, i \in [{3k+2  \over 2}]_0,  j\in [3]_0 \}
\cup \{( (i  + {3k+2 \over 2} +  j k) \, {\rm mod} \, (3k+1)), 2i+1) \, | \, i \in [{3k \over 2}]_0,  j\in [3]_0 \}$. 
For $k=6$ the construction is depicted on the right-hand side of Fig.  \ref{C17C17}. 
We can see that
$M = \{(i,  (2i + j(k+1)) \, {\rm mod} \, (3k+1)) \, | \, i \in [3k+1]_0,  j\in [3]_0 \} $.

Note that the distance between vertices of $M \cap C_{3k+1}^i $ is either $k$ or  $k+1$, while the distance 
between vertices of  $M \cap ^{j}\!C_{3k+2}$ is either $k+1$ or  $k-1$.
Thus,  using the same arguments as in the proof of Proposition \ref{mod3}  we can see that $M$ satisfies conditions (i), (ii), and (iii) of Lemma \ref{l1}.

Since we can confirm that $d(u,v) \ge 3$ for every $u,v \in M$,  
It follows that $M $ is a mutual-visibility set of cardinality  $9k+3$ in $C_{ 3k+1}  \Box C_{3k+1}$ if $k\ge 6$ is even.
\end{proof}

\begin{prop} \label{mod2}
Let $t\equiv 2 $ (mod $6$).  If $t \ge 20$,  then $\mu(C_t \Box C_t) = 3t$.
\end{prop}

\begin{figure}[H]
     \centering
	\includegraphics[width=12cm]{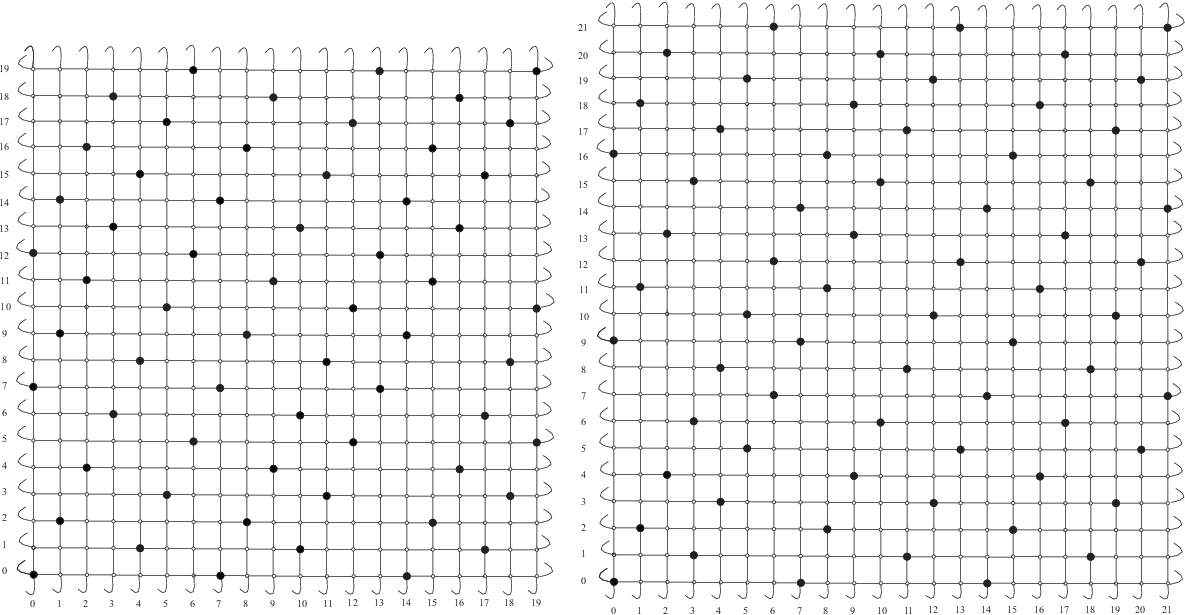}
	\bigskip
	\caption{Mutual-visibility set of $C_{20} \Box  C_{20}$  (left) and   $C_{22} \Box  C_{22}$  (right)}
\label{C20C20}
\end{figure}

\begin{proof}
We will show that a mutual-visibility set of cardinality  $3t$ exists for every $t \equiv 2$ (mod $6)$, $t \ge 20$.  

Let us of an even $k$ set $t=3k+2$.
Remind that for every $u,v \in V(C_{3k+2} \Box C_{3k+2})$  we have 
$I(u,v) = [u_x,v_x]_{3k+2} \times [u_y,v_y]_{3k+2}$. Moreover,  the graph induced by
$[u_x,v_x]_{3k+2} \times [u_y,v_y]_{3k+2}$  is isomorphic to a subgraph of 
$P_{ {3k+4 \over 2}}  \Box P_{3k+4 \over 2}$.

Let $M =  \{((i  + j (k+1)) \, {\rm mod} \, (3k+2), 2i) \, | \, i \in [{3k+2  \over 2}]_0,  j\in [3]_0 \}
\cup \{ ((i  + {3k +2 \over 2} +  j (k+1)) \, {\rm mod} \, (3k+2), 2i+1)) \, | \, i \in [{3k+2 \over 2}]_0,  j\in [3]_0 \}$.
For $k=6$, the construction is depicted on the left-hand side of Fig.  \ref{C20C20}.

Since
 $M = \{( i, 2i ) \, | \, i \in [k+1]_0 \} \cup \{(i, 2i + k+1) \, | \, i \in [k+1]_0 \}
\cup \{(i, (2i + 2k)  \, {\rm mod} \, (3k+1)) \, | \, i \in [k+1]_0 \} 
\cup \{( i+k+1, 2i) \, | \, i \in [k+1]_0 \} \cup \{(i +k+1, 2i + k +1) \, | \, i \in [k+1]_0 \}
\cup \{( i + k +1, (2i + 2k + 2)  \, {\rm mod} \, (3k+1)) \, | \, i \in [k+1]_0 \} 
\cup \{(i + 2k+2, 2i) \, | \, i \in [k]_0 \} \cup \{(i +2k+2, 2i + k+3  ) \, | \, i \in [k]_0 \}
\cup \{(i + 2k+2, (2i + 2k+2)  \, {\rm mod} \, (3k+1)) \, | \, i \in [k]_0 \}  $,
the vertices of 
$M \cap ^{j}\!C_{3k+2}$  are pairwise  at distance  from the set $\{ k+3,  k, k-1\}$,
while the distance between vertices of  
$M \cap C_{3k+2}^i $ is  either $k+1$ or  $k$.
Thus,  using the same arguments as in the proof of Proposition \ref{mod3}  we can see that $M$ satisfies conditions (i), (ii), and (iii) of Lemma \ref{l1}.

Let  $A =\{(0,0),  (k+1,0),  (2k+2, 0) \}$ and $B =\{ (3k+1,3k+1),  (k, 3k+1),  (2k+1,3k+1) \}$.
We can see that $A \cup B \subseteq M$.  Moreover, 
$d((0,0),  (3k+1,3k+1)) = d((k+1,0),  (k,3k+1)) = d((2k+2,0),  (2k+1,3k+1))=2$, while for every 
$u,v \in M \setminus A $ as well as $u,v \in M \setminus B$ we have $d(u,v)\ge 3$.  Thus,   $M \setminus A$ and
$M \setminus B$ are both mutual-visibility sets in $C_{3k+2} \Box C_{3k+2}$. 

If $x \in A$ and  $y \in B$  such that $d(x,y)=2$, then we say that $x$ and $y$ form a {\em cross pair} in $M$. 
We have to show that for every $u, v \in M$, where $I(u,v)$ admits a cross pair, there exists an $M$-unhindered  
$u,v$-path in  $C_t \Box C_t$. 

 Suppose w.l.o.g. that $u_y > {3k+2  \over 2}$ and $v_y < {3k+2  \over 2}$. 
Assume first that  $u_x \le 2k+1$. 
Let   $x=(2k+1,3k+1)$, $y= (2k+2, 0)$  and $x,y  \in I(u,v)$.  Let also $x'=(2k,3k+1)$, $x''=(2k+1,3k)$, 
$y'= (2k+3, 0)$ and $y''= (2k+2, 1)$.  Note that every $M$-unhindered  $u,x$-path contains either $x'$ or $x''$,
while every $M$-unhindered   $v,y$-path contains either $y'$ or $y''$.  

If $u_y  \not \in \{3k-1,3k\}$ and $u_x < 2k$, it is not difficult to see that we can always find an $M$-unhindered 
shortest $u,x'$- and $u,x''$-path. Analogously, if  $v_y > 1$ and $v_x \not \in \{2k+2, 2k+3k\}$, we can always find 
an $M$-unhindered   $v,(2k+2, 1)$- and $v,(2k+3, 0)$-path.  It follows that we can construct an 
an $M$-unhindered  $u,v$-path for every $u,v$ that satisfy the above conditions.

If  $u_y = 3k + 1$ or $u_y = 3k$ and $u_x =  {3 k \over 2}$,  then $v_y > 0$ since the vertices of 
$M \cap C_{3k+2}^i $ are pairwise  at distance either $k+1$ or  $k$. 
Thus,  there exists an $M$-unhindered  $u,v$-path that contains  $(2k+2, 1)$.
Analogously, if $v_y = 0$ or $v_y = 1$ and $v_x =  {k \over 2} + 1$, then $u_y < 3k+1$ 
since the vertices of  $M \cap ^{j}\!C_{3k+2}$ are pairwise at distance from the set $\{ k+3,  k, k-1\}$.
Thus, there exists an $M$-unhindered  $u,v$-path that contains  $(2k+1,3k)$.

Finally, if $u_y = 3k$ and $u_x =  {3 k \over 2}$ then we can an find an $M$-unhindered 
 $u,x'$- and $u,x''$-path, while for 
$v_y =  1$ and $v_x = {5k \over 2}+3$ we can  find an $M$-unhindered  $v,y'$- and $v,v''$-path. 
It follows that we can construct an $M$-unhindered  $u,v$-path for this selection of $u$ and $v$.

The above discussion shows that there exists an $M$-unhindered  $u,v$-path for every $u,v \in M$ 
where $u_x \le 2k+1$.  The proof for  $u_x \ge 2k+2$ is analogous. 

Since   for every $u, v \in M$, where  
either $\{ (0,0),  (3k+1,3k+1) \} \subseteq I(u,v)$ or $\{ (k+1,0),  (k,3k+1)  \} \subseteq I(u,v)$,
we can establish in a similar fashion that there exists an $M$-unhindered  
$u,v$-path in  $C_t \Box C_t$,  the proof is complete.   
\end{proof}

\begin{prop} \label{mod4}
Let $t\equiv 4 $ (mod $6$).  If $t \ge 22$,  then $\mu(C_t \Box C_t) = 3t$.
\end{prop}

\begin{proof}
Since  $\mu(C_t \Box C_t) \le 3t$ by Proposition \ref{p3},  we have to show that a mutual-visibility set of cardinality  $3t$ exists for every $t \equiv 4$ (mod $6)$.  

Let us of an odd $k$ set $t=3k+1$.
Let 
$M =  \{((i  + j k) \, {\rm mod} \, (3k+1)),  2i) \, | \, i \in [{3k+1  \over 2}]_0,  j\in [3]_0 \}
\cup \{((i  + {3k+1 \over 2} +  j k) \, {\rm mod} \, (3k+1)), 2i+1) \, | \, i \in [{3k+1 \over 2}]_0,  j\in [3]_0 \}$. 
For $k=6$ the construction is depicted on the right-hand side of Fig.  \ref{C20C20}. 

We can see that
$M = \{(i, 2i) \, | \, i \in [k]_0 \} \cup \{(i, 2i + k +2) \, | \, i \in [k]_0 \}
\cup \{(i, (2i + 2k +2)  \, {\rm mod} \, 3k ) \, | \, i \in [k]_0 \} 
\cup \{(i+k,  2i) \, | \, i \in [k]_0 \} \cup \{(i +k,  2i + k +2) \, | \, i \in [k]_0 \}
\cup \{(i + k, (2i + 2k)  \, {\rm mod} \, 3k) \, | \, i \in [k]_0 \} 
\cup \{(i + 2k, 2i) \, | \, i \in [k+1]_0 \} \cup \{( i +2k, 2i + k) \, | \, i \in [k+1]_0 \}
\cup \{( i + 2k, (2i + 2k)  \, {\rm mod} \, 3k) \, | \, i \in [k+1]_0 \}  $.
It follows that the distance between vertices 
of  $M \cap ^{j}\!C_{3k+2}$ is  from the set $\{ k-1,  k,  k+2\}$,  while
the vertices of $M \cap C_{3k+1}^i $ are pairwise  at distance either $k+1$ or  $k$.
Thus,  using the same arguments as in the proof of Proposition \ref{mod3}  we can establish that $M$ satisfies conditions (i), (ii), and (iii) of Lemma \ref{l1}.

Note that $M$ admits three cross pairs:  $\{(0,0),  (3k+1,3k+1)\}$,  $\{(k,0),  (k-1,3k+1)\}$ and
$\{(2k,0),  (2k-1,3k+1)\}$.  Analogously as in the proof of Proposition \ref{mod2}, we can now prove that for every
$u,v \in M$ there exists an $M$-unhindered  $u,v$-path in  $C_t \Box C_t$. 
\end{proof}

\begin{thm} \label{t2}
Let $t \ge 14$ or $t=12$.  If $s \ge t$,  then $\mu(C_s \Box C_t) = 3t$.
\end{thm}

\begin{proof}
Let $t \equiv 3$ (mod $6$).
If $s=t$,  this is Proposition \ref{p4}.  Since the upper bound is given by Proposition \ref{p4}, 
we have to show that there exists a mutual-visibility set of cardinality  $3t$ for every $s>t$.

Let $M = \{((2i + jk) \, {\rm mod} \, 3k, i) \, | \, i \in [3k]_0,  j\in [3]_0 \}$.
Remind that $M$ is a mutual-visibility set of cardinality  $3t$ in $C_t \Box C_t$ if $t \equiv 3$ (mod $6$). 

For $i  \in [3k]$,  let  $S_i = \{ (j \, {\rm mod} \, 3) k + (\lfloor { j \over 3} \rfloor   \, {\rm mod} \, k)  \, | \, 0 \le j < i \}$.
Note that  $S_i \subseteq [3k]_0$.  Moreover, elements of  $S_i$ are uniformly distributed among sets 
$[k]_0$,  $[k,2k-1]$ and $[2k,3k-1]$, i.e., 
each of the sets   $[k]_0$,  $[k,2k-1]$ and $[2k,3k-1]$ contains 
 either  $\lfloor { i \over 3} \rfloor$ or  $\lfloor { i \over 3} \rfloor + 1 $ elements from $S_i$. 
For a set $S \subseteq [3k]_0$, let the function $\delta_S :  [3k]_0 \rightarrow [3k]_0$ be defined by $\delta_S(i) =  |\{ j  \, | \,  j \in S \, {\rm and } \, j < i  \}|$.  
Moreover,  let $M_i = \{ (u_x + \delta_{S_i}(u_x), u_y)  \, | \,  u \in M \}$.  
In other words,  a vertex of $M_i$ is obtained from $u \in M$ by shifting its first coordinate by $\delta_{S_i}(u_x)$ positions.  Note that $M_i$ is a subset of $C_{t+i} \Box C_{t}$ of cardinality  $3t$.  

Loosely speaking,  we can say that $M_i$ is obtained from $M$ by inserting $i$ empty $C_{3k}$-fibers.  
To show that $M_i$ is a mutual-visibility set  in $C_{t+i} \Box C_{t}$, first note that   condition  (iv) of Lemma \ref{l1}
clearly holds.  For conditions (i)-(iii) we apply the fact that  elements of $S_i$ are uniformly distributed among sets 
$[k]_0$,  $[k,2k-1]$ and $[2k,3k-1]$.  Since for every $j \in [k]_0$ and every $u,v \in M \cap ^jC_{s}$ we 
have $d_{C_t \Box C_{t}}(u,v) = k$,  it follows that for $u,v \in M_i \cap ^jC_{s+i}$ it holds that
either $d_{C_{t+i} \Box C_{t}}(u,v) = k+\lfloor { i \over 3} \rfloor $ or 
$d_{C_{t+i} \Box C_{t}}(u,v) = k+\lfloor { i \over 3} \rfloor + 1$.  Since for every 
$z,w \in V(C_{ 3k}  \Box C_{3k})$ we have $d(z_x,w_x)  \le {3k + i  \over 2}$,   
conditions (i)-(iii) of Lemma \ref{l1} are also fullfiled. Thus, $M_i$   is a mutual-visibility set  in $C_{t+i} \Box C_{t}$ for 
 $i  \in [3k]$. In order to obtain a mutual-visibility set  in $C_{t+i} \Box C_{t}$ for $i > 3k$,  we consecutively 
 apply the above procedure for $M=M_{2t}$, $M=M_{3t}, \ldots$  
 
For other five cases, i.e.  if $t \equiv r$ (mod $6$), where  $r \in  \{0, 1, 2, 4, 5 \}$, the proof is analogous. We start 
with the mutual-visibility set  $M$ of cardinality  $3t$ of $C_t \Box C_t$ where  $t \equiv r$ (mod $6$).  
Next,  the sets $S_i$ are defined for every $i \in [t]$. Finally,  the sets $M_i$ are determined for every $i \in [t]$.
In particular,
\begin{itemize}
\item
 for $t \equiv 0$ (mod $6$),  we set $M =  \{((i  + j k, 2i + \ell k) \, | \, i \in [{k  \over 2}]_0,  j, \ell \in [3]_0 \}
\cup \{(k-i-1+jk,  (2i  + 3 + \ell k ) \, {\rm mod} \, 3k) \, | \, i \in [{k \over 2}]_0,  j, \ell \in [3]_0 \}$  and
for $i  \in [3k]$,  we set  $S_i = \{ (j \, {\rm mod} \, 3) k + (\lfloor { j \over 3} \rfloor   \, {\rm mod} \, k)  \, | \, 0 \le j < i \}$;

\item
 for $t \equiv 5$ (mod $6$),  then we set 
$M =  \{(i  + j (k+1) \, {\rm mod} \, (3k+2)),  2i ) \, | \, i \in [{3k+3  \over 2}]_0,  j\in [3]_0 \}
\cup \{(( i  + {3k+3 \over 2} +  j (k+1) \, {\rm mod} \, (3k+2)), 2i+1) \, | \, i \in [{3k+1 \over 2}]_0,  j\in [3]_0 \}$ and 
for $i  \in [3k+2]$ we set  $S_i = \{ (j \, {\rm mod} \, 3) (k+1) + (\lfloor { j \over 3} \rfloor   \, {\rm mod} \, (k+1) )  \, | \, 0 \le j < i \}$;

\item
 for $t \equiv 1$ (mod $6$),  let 
$M =  \{((i  + j k) \, {\rm mod} \, (3k+1)),2i) \, | \, i \in [{3k+2  \over 2}]_0,  j\in [3]_0 \}
\cup \{( (i  + {3k+2 \over 2} +  j k) \, {\rm mod} \, (3k+1)), 2i+1) \, | \, i \in [{3k \over 2}]_0,  j\in [3]_0 \}$ and 
 for $i  \in [3k]$,  let  $S_i = \{ (j \, {\rm mod} \, 3) k + (\lfloor { j \over 3} \rfloor   \, {\rm mod} \, k)  \, | \, 0 \le j < i \}$
 and $S_{3k+1} =[3k+1]_0$;

\item
 for $t \equiv 2$ (mod $6$),  let
$M =  \{((i  + j (k+1)) \, {\rm mod} \, (3k+2), 2i) \, | \, i \in [{3k+2  \over 2}]_0,  j\in [3]_0 \}
\cup \{ ((i  + {3k +2 \over 2} +  j (k+1)) \, {\rm mod} \, (3k+2), 2i+1)) \, | \, i \in [{3k+2 \over 2}]_0,  j\in [3]_0 \}$ and 
for $i  \in [3k+2]$,  let  $S_i = \{ (j \, {\rm mod} \, 3) (k+1) + (\lfloor { j \over 3} \rfloor   \, {\rm mod} \, (k+1) )  \, | \, 0 \le j < i \}$;

\item
 for $t \equiv 4$ (mod $6$),  let 
$M =  \{((i  + j k) \, {\rm mod} \, (3k+1)),  2i) \, | \, i \in [{3k+1  \over 2}]_0,  j\in [3]_0 \}
\cup \{((i  + {3k+1 \over 2} +  j k) \, {\rm mod} \, (3k+1)), 2i+1) \, | \, i \in [{3k+1 \over 2}]_0,  j\in [3]_0 \}$
and  for $i  \in [3k]$,  let  $S_i = \{ (j \, {\rm mod} \, 3) k + (\lfloor { j \over 3} \rfloor   \, {\rm mod} \, k)  \, | \, 0 \le j < i \}$
 and $S_{3k+1} =[3k+1]_0$;
\end{itemize}

We are left to show that $C_s \Box C_t$  admits a mutual-visibility set of cardinality $3t$,  if 
$t \in  \{12, 14,  16 \}$ and $s \ge t$.  

If  $t=12$,  we found a mutual-visibility set of cardinality $3t$ for every $s \in  \{12, 13,  14,  15,  16,  17,  18 \}$. 
Moreover,  a mutual-visibility set for $C_{18} \Box C_{12}$ can be   applied in $C_{s} \Box C_{12}$ for 
every $s > 18$ analogously as in the proof for $t \equiv 3$ (mod $6$).

If  $t=14$,  we found a mutual-visibility set of cardinality $3t$ for every $s \in  \{14,  15,  16,  17,  18, 19, 20 \}$. 
Moreover,  a mutual-visibility set for $C_{20} \Box C_{14}$ can be   applied in $C_{s} \Box C_{14}$ for 
every $s > 20$ analogously as in the proof for $t \equiv 2$ (mod $6$).

If  $t=16$,  we found a mutual-visibility set of cardinality $3t$ for every $s \in  \{16,  17,  18 \}$. 
Moreover,  a mutual-visibility set for $C_{18} \Box C_{16}$ can be   applied in $C_{s} \Box C_{16}$ for 
every $s > 20$ analogously as in the proof for $t \equiv 3$ (mod $6$).

The above-mentioned constructions  are available on the web page 
\href{https://omr.fnm.um.si/wp-content/uploads/2023/08/MVS-constructions-for-two-cycles.txt}{external document}.

Since the presented solutions comprise all graphs of interest,  the proof is complete.
\end{proof}

\begin{table}[H] 
\begin{center}
\caption{Mutual-visibility numbers of $C_s \Box C_t$,  where $t \le 13 $}  \label{tab2}
\bgroup
\def\arraystretch{1.5}
\begin{tabular}{|c||c|c|c|c|c|c|c|c|c|c|c|c|}
  \hline
  $t$ \textbackslash $s$ & $3$ & $4$ & $5$ & $6$ & $7$  &  $8$  &  $9$ &  $10$ &  $11$ &  $12$ &  $13$ &  $14$   \\ 
  \hline
 \hline
  $3$ & 6  & 7  & 7  & 9  & \ditt & \ditt & \ditt & \ditt  & \ditt  & \ditt &  \ditt & \ditt    \\
\hline
  $4$ &  &  $9$  & $10$  & $11$  & $11$  & $12$  & \ditt & \ditt & \ditt  & \ditt&  \ditt & \ditt  \\
\hline
  $5$ & & &  10 & 12 & 13 & 15 & \ditt & \ditt & \ditt &  \ditt & \ditt &  \ditt   \\
\hline
 $6$ &  &  &  & $14$  & $15$ & $17$  & $18$  & \ditt & \ditt & \ditt &  \ditt & \ditt  \\
\hline
  $7$ & &  &  &  & $16$  & $18$  & $18$  & $20$  & $20$  & $21$ &  \ditt & \ditt  \\
\hline
  $8$ & & &  &  &  & $21$  & $21$  & $23$  & $23$ &  $24$ &  \ditt & \ditt   \\
\hline
  $9$ &  & & &  &  &  & $22$  & $25$  & $25$  & $27$ & \ditt & \ditt   \\
    \hline
  $10$ &  & & &  &  &  & &  $27$  & $27$  & $30$   &  \ditt  &  \ditt   \\
    \hline
  $11$ & &  & & &  &  &  & &  $29$  & $32$  & $33$ &  \ditt  \\
    \hline
  $12$ & & &  & & &  &  &  & & $36$  &  \ditt &  \ditt  \\
    \hline
  $13$ & & & &  & & &  &  &  & &  $38$  & $39$   \\
    \hline
\end{tabular}
\egroup
\end{center}
\end{table}

The next proposition shows that  even for smaller $t$  it holds that 
$\mu(C_s \Box C_t) = 3t$ providing that $s$ is big enough. 

\begin{prop} \label{CsCtMali}
Let $(s_r, t) \in \{ (6, 3), (8,4), (8,5), (9,6),  (12,7),  (12,8),  (12,9), $ $(12,10),   (13,11),   (14,  13)  \}$.   If
$s \ge s_r$,  then $\mu(C_s \Box C_t) = 3t$.
\end{prop}

\begin{proof}  
The proof is based on the constructions of mutual-visibility sets for the graphs of interest that are 
available on the web page \href{https://omr.fnm.um.si/wp-content/uploads/2023/08/MVS-constructions-for-two-cycles.txt}{external document}.

We do not give the details here since the reasoning is analogous for all values of $t$.
As an example,  we consider the case $t=3$.  We found a  mutual-visibility set with 9 vertices of $C_s \Box C_3$ for every 
$s \in [6,18]$.  Moreover, we can apply the construction for $C_{18} \Box C_3$ in a similar manner to the 
	case $t \equiv 3$ (mod $6$) in the proof of Theorem \ref{t1}. 
\end{proof}

\begin{rem}  
Mutual-visibility numbers of $C_s \Box C_t$,  for $t \le 13 $ are presented in 
Table   \ref{tab2}.   
The majority of the results were obtained using a computer, especially through backtracking.
Note that the last mutual-visibility number in a row that corresponds to some $t$  always equals $3t$.
\end{rem}

\section*{Data availability}
The datasets generated during and/or analysed during the current study are available on the web pages 
\href{https://omr.fnm.um.si/wp-content/uploads/2023/08/MVS-constructions-for-paths-and-cycles.txt}{external document} and 
\href{ https://omr.fnm.um.si/wp-content/uploads/2023/08/MVS-constructions-for-two-cycles.txt}{external document}
as well as from the corresponding author on reasonable request.

\section*{Acknowledgments}
This work was supported by the Slovenian Research Agency under the grants P1-0297,  J1-2452, and  J2-1731.


\end{document}